\documentclass[twoside,11pt,reqno]{article}
\usepackage{amssymb}
\usepackage{amsthm}
\usepackage[tbtags]{amsmath}
\usepackage{doc}
\usepackage{latexsym}
\usepackage{amscd}
\font\headd=cmr8  \theoremstyle{change}

  \font\headd=cmr8
\theoremstyle{change}  
\begin{document} \thispagestyle{plain}
 \markboth{}{}
\small{\addtocounter{page}{0} \pagestyle{plain} \vspace{0.2in}
\pagestyle{myheadings}
 \markboth{\headd Taekyun Kim $~~~~~~~~~~~~~~~~~~~~~~~~~~~~~~~~~~~~~~~~~~~~~$}
 {\headd $~~~~~~~~~~~~~~~~~~~~~~~~~~~~~~~~~~~~~~~~~~$Euler numbers and polynomials}
\noindent{\large \bf $p$-adic $l$-functions and sums of powers}
\footnote{{}
}
\vspace{0.15in}\\
\noindent{\sc Taekyun Kim}
\newline
{\it Jangjeon Research Institute for Mathematical Sciences $\&$
Physics, Ju-Kong Building 103-Dong 1001-ho, 544-4 Young-chang Ri
Hapcheon-Up Hapcheon-Gun Kyungnam, 678-802, Korea\\
e-mail} : {\verb|tkim64@hanmail.net|} ~{\it or}~
{\verb|tkim@kongju.ac.kr|}
\vspace{0.15in}\\
{\footnotesize {\sc Abstract.} In this paper, we give an explicit
$p$-adic expansion of
$$\sum_{\substack{j=1\\ (j, p)=1}}^{np} \dfrac{(-1)^{j}}{j^{r}}$$
as a power series in $n$. The coefficients are values of $p$-adic
$l$-function for Euler numbers.}
\vspace{0.2in}\\
\noindent{\bf 1. Introduction} \setcounter{equation}{0}
\vspace{0.1in}\\
\indent Let $p$ be a fixed prime. Throughout this paper
$\mathbb{Z}_{p}$, $\mathbb{Q}_{p}$, $\mathbb{C}$ and
$\mathbb{C}_{p}$ will, respectively, denote the ring of $p$-adic
rational integers, the field of $p$-adic rational numbers, the
complex number field and the completion of algebraic closure of
$\mathbb{Q}_{p}$, cf. [1], [3], [6], [10]. Let $v_{p}$ be the
normalized exponential valuation of $\mathbb{C}_{p}$ with
$|p|_{p}=p^{-v_{p}(p)}= p^{-1}$. Kubota and Leopoldt proved the
existence of meromorphic functions, $L_{p}(s, \chi)$, defined over
the $p$-adic number field, that serve as $p$-adic equivalents of
the Dirichlet $L$-series, cf. [8], [10]. These $p$-adic
$L$-functions interpolate the values
$$L_{p}(1-n, \chi) = -\dfrac{1}{n}(1-\chi_{n}(p) p^{n-1}) B_{n, \chi_{n}}, \quad {\rm
for}~ n \in \mathbb{N}=\{1, 2, \cdots\}, $$ where $B_{n, \chi}$
denote the $n$th generalized Bernoulli numbers associated with the
primitive Dirichlet character $\chi$, and $\chi_{n} = \chi w^{-n}$,
with $w$ the {\it Teichm\"uller} character, see [2, 3, 5, 6, 17,
20]. In [14], L. C. Washington have proved the below interesting
formula:
$$\sum_{\substack{j=1 \\ (j, p)=1}}^{np} \dfrac{1}{j^{r}}=-\sum_{k=1}^{\infty} \binom{-r}{k} (pn)^{k} L_{p}(r+k, ~w^{1-k-r}),$$
where $\binom{-r}{k}$ is binomial coefficient. In the recent many
authors have studied $q$-extension of Euler numbers and Bernoulli
numbers (see [1, 4, 5, 9, 12, 13]). These $q$-extensions seem to be
valuable and worthwhile in the areas of mathematical physics and
mathematics (see [ 1, 4, 6, 13, 14, 15, 17, 18, 19]). By using
$q$-Volkenborn integration, Kim gave the interesting properties of
$q$-Bernoulli and Euler polynomials [8, 9, 10, 11] and
Ryoo-Kim-Agarwal have investigated the properties of the
$q$-extension of Euler numbers and polynomials by using
``Mathematica package", see [14, 15]. The problems to find the sums
of powers of consecutive $q$-integers were suggested. Kim and
Schlosser treated the formulae for the sums of powers of consecutive
$q$-integers [7,11, 16] and these formulae were used to give the
$q$-extension of Washington's $p$-adic $L$-functions and sums of
powers ( see [ 6, 10, 20]). In [11, 14], we found the interesting
formulae ``alternating sums of powers of consecutive integers "
which are related to Euler numbers and polynomials. By using these
alternating sums of powers of consecutive integers, we try to
construct the $p$-adic $l$-functions and sums of powers for Euler
numbers and polynomials, corresponding to Washington and Kim (see
[10, 20]). The purpose of this paper is to give alternating $p$-adic
harmonic series in terms of $n$ and $p$-adic $l$-function for Euler
numbers.
\vspace{0.2in}\\
\noindent{\bf 2. A note on $l$-series associated with Euler
numbers and polynomials}
\vspace{0.1in}\\
\indent  We begin with well known Euler polynomials $E_{n}(x)$.
\vspace{0.1in}\\
\noindent {\bf Definition 1.} Euler polynomials are defined by
$$\dfrac{2}{e^{t}+1} e^{xt} = \sum_{n=1}^{\infty} \dfrac{E_{n}(x)}{n!} t^{n},$$
$E_{n}(x)$ are called $n$-th Euler polynomials. For $x=0$,
$E_{n}=E_{n}(0)$ are called Euler numbers. By the definition of
Euler polynomials, we easily see that,
$$E_{l}(x) = \sum_{n=0}^{l} \binom{l}{n} E_{n} x^{l-n} ~\in
\mathbb{C}[x].$$ From the generating function of Euler polynomials
$F(t, x) = \frac{2}{e^{t}+1} e^{xt}$, we derive
\begin{eqnarray}
F(t, x) &=& 2e^{xt} \sum_{l=0}^{\infty} (-1)^{l} e^{lt} \notag\\
&=& 2 \sum_{l=0}^{\infty} (-1)^{l} e^{(l+x)} t.
\end{eqnarray}

For $k\in \mathbb{N}$, we note that
\begin{eqnarray}
\dfrac{d^{k}}{dt^{k}} F(t, x) \Big|_{t=0} &=& 2
\sum_{l=0}^{\infty}
(-1)^{l} \dfrac{d^{k}}{dt^{k}} e^{(l+x)t} \Big|_{t=0}\notag\\
&=& 2 \sum_{l=0}^{\infty} (-1)^{l} (l+x)^{k}.
\end{eqnarray}
Therefore we can define the Euler zeta function as follows:
\vspace{0.1in}\\
\noindent {\bf Definition 2.} For $s \in \mathbb{C}$, we define
Euler zeta function as
\begin{equation}
\zeta_{E}(s) = 2 \sum_{l=0}^{\infty} \dfrac{(-1)^{l}}{(l+x)^{s}}.
\end{equation}
By using Definition 2 and (2), we obtain the following:
\vspace{0.1in}\\
\noindent {\bf Proposition 3.} {\it For $k \in \mathbb{N}$, we
have}
\begin{equation}
\zeta_{E}(-k, x) = E_{k}(x).
\end{equation}

For $f$(=odd) $\in \mathbb{N}$,
\begin{eqnarray*}
\sum_{n=0}^{\infty} E_{n}(x) \dfrac{t^{n}}{n!} &=&
\dfrac{2}{e^{t}+1} e^{xt} = 2 \sum_{l=0}^{\infty} (-1)^{l}
e^{(l+x)t} = 2 \sum_{a=0}^{f-1} \sum_{l=0}^{\infty} (-1)^{a+lf}
e^{(a+lf+x)t}\\
&=& \sum_{a=1}^{f} (-1)^{a} 2\sum_{n=0}^{\infty} (-1)^{l}
e^{ft(l+\frac{x+a}{f})} = \sum_{a=1}^{f} (-1)^{a}
\sum_{n=0}^{\infty} E_{n} \left(\dfrac{x+a}{f}\right)
\dfrac{f^{n}t^{n}}{n!}\\
&=& \sum_{n=0}^{\infty} f^{n} \sum_{a=1}^{f} (-1)^{a} E_{n}
\left(\dfrac{x+a}{f} \right) \dfrac{t^{n}}{n!}.
\end{eqnarray*}
Thus we note that
\begin{equation}
E_{n}(x) = f^{n} \sum_{a=1}^{f} (-1)^{n} E_{n}\left(\dfrac{x+a}{f}
\right), \end{equation} where $f$(=odd) $\in \mathbb{N}$. This (5)
is so called {\it Distribution for Euler polynomials}.
$$2\sum_{l=0}^{n-1}(-1)^{l} l^{m} = (-1)^{n+1} \sum_{l=0}^{m-1}
E_{l} n^{m-l}\binom{m}{l}+\left( (-1)^{n+1}+1 \right) E_{m}.$$

In particular, if $n$ is even, then
$$2\sum_{l=0}^{n-1}(-1)^{l}
l^{m} = - \sum_{l=0}^{m-1} E_{l} n^{m-l}.$$

Let $s$ be a complex variable and let $a,~F$(=odd) be integers
with $0 < a < F$.
\begin{eqnarray}
H(s,~ a|F) &=& \sum_{\substack {m\equiv a(F) \\ m > 0}}
\dfrac{(-1)^{m}}{m^{s}} = \sum_{n=0}^{\infty}
\dfrac{(-1)^{nF+a}}{(a+nF)^{s}} = (-1)^{a} \sum_{n=0}^{\infty}
\dfrac{(-1)^{n}}{(a+nF)^{s}}\notag\\
& =& (-1)^{a} \sum_{n=0}^{\infty}
\dfrac{(-1)^{n}}{F^{s}\left(\dfrac{a}{F}+n\right)^{s}} =
\dfrac{(-1)^{a} F^{-s}}{2}  ~2\sum_{n=1}^{\infty}
\dfrac{(-1)^{n}}{n+\dfrac{a}{F}} \notag\\
&=& \dfrac{(-1)^{a}F^{-s}}{2} ~\zeta_{E} \left(s,
\dfrac{a}{F}\right).
\end{eqnarray}

Note that
\begin{equation}
H(-n,~a|F)= (-1)^{a} \dfrac{F^{n}}{2} E_{n} \left( \dfrac{a}{F}
\right).
\end{equation}

Let $\chi$ be the primitive Dirichlet character with conductor
$f$(=odd) $\in \mathbb{N}$.
\begin{eqnarray}
F_{\chi}(t) &=& 2\sum_{n=0}^{\infty} e^{nt} \chi(n) (-1)^{n}\notag\\
&=& 2\sum_{a=0}^{f-1} \sum_{n=0}^{\infty} e^{(a+nf)t} \chi(a+nf) (-1)^{a+nf} \notag\\
&=& \sum_{a=0}^{f-1} e^{at} \chi(a) (-1)^{a} 2 \sum_{n=0}^{\infty} (-1)^{n}e^{nft}\notag\\
&=& 2 \sum_{a=0}^{f-1} e^{at} \chi(a) (-1)^{a} \left(
\dfrac{1}{e^{ft}+1}\right)\notag\\
&=& \dfrac{2 \sum_{a=0}^{f-1} e^{at} \chi(a) (-1)^{a}}{e^{ft}+1}\notag\\
&=& \sum_{n=0}^{\infty} E_{n, ~\chi}~\dfrac{t^{n}}{n!}.
\end{eqnarray}
Thus, we can define the below generalized Euler number attached to
$\chi$.
\vspace{0.1in}\\
\noindent {\bf Definition 4.} Let $\chi$ be the Dirichlet
character with conductor $f$(=odd) $\in \mathbb{N}$. Then we
define the generalized Euler numbers attached to $\chi$ as
follows;
$$\dfrac{2 \sum_{a=0}^{f-1} e^{at} \chi(a) (-1)^{a}}{e^{ft}+1} =
\sum_{n=0}^{\infty} E_{n, \chi}~\dfrac{t^{n}}{n!}.$$ $E_{n, \chi}$
will be called the $n$-th generalized Euler numbers attach to
$\chi$.

\medskip

From the Definition 4, we derive the below formula:
\begin{eqnarray*}
\sum_{n=0}^{\infty} E_{n, \chi}~\dfrac{t^{n}}{n!} &=& \dfrac{2
\sum_{a=0}^{f-1} e^{at} \chi(a) (-1)^{a}}{e^{ft}+1} =
\sum_{a=0}^{f-1} \chi(a) (-1)^{a} \left( \dfrac{2}{e^{ft}+1}
e^{at} \right)\\
&=& \sum_{a=0}^{f-1} \chi(a) (-1)^{a} \sum_{n=0}^{\infty}
E_{n}\left(\dfrac{a}{f}\right) \dfrac{f^{n}t^{n}}{n!} =
\sum_{n=0}^{\infty} f^{n} \sum_{a=0}^{f-1} \chi(a) (-1)^{a} E_{n}
\left(\dfrac{a}{f} \right) \dfrac{t^{n}}{n!}.
\end{eqnarray*}
By comparing the coefficients on both sides, we easily see that
\begin{equation}
E_{n,\chi} = f^{n} \sum_{a=0}^{f-1} \chi(a) (-1)^{a}
E_{n}\left(\dfrac{a}{f} \right). \end{equation}
\vspace{0.1in}\\
\noindent {\bf Definition 5.} For $s \in \mathbb{C}$, we define
Dirichlet's $l$-function as follows: $$l(s, ~\chi) = 2
\sum_{n=1}^{\infty} \dfrac{\chi(n) (-1)^{n}}{n^{s}}.$$

Note that
\begin{eqnarray*}
F_{\chi}(t) &=& 2 \sum_{n=1}^{\infty} e^{nt}\chi(n) (-1)^{n} =
\dfrac{2 \sum_{a=0}^{f-1} (-1)^{a}\chi(a)e^{at}}{e^{ft}+1} =
\sum_{n=0}^{\infty} E_{n, \chi} ~\dfrac{t^{n}}{n!}.
\end{eqnarray*}

For $k \in \mathbb{N}$,
\begin{eqnarray}
E_{k, \chi} &=& \dfrac{d^{k}}{dt^{k}} F_{\chi}(t) \Big|_{t=0} = 2
\sum_{n=1}^{\infty} \chi(n) (-1)^{n} \dfrac{d^{k}}{dt^{k}}
e^{nt} \Big|_{t=0}\notag\\
&=& 2 \sum_{n=1}^{\infty} \chi(n) (-1)^{n} n^{k}.
\end{eqnarray}
By Definition 5 and (10), we easily see that $l(-k, ~\chi)= E_{k,
~\chi}$, where $k \in \mathbb{N}$. Therefore we obtain the
following:
\vspace{0.1in}\\
\noindent {\bf Proposition 6.} {\it Let $k$ be the positive
integer. Then we have}
\begin{equation}
l(-k, \chi) = E_{k, \chi}.
\end{equation}

\medskip

Let $\chi$ be the Dirichlet character with conductor $f$(=odd)
$\in \mathbb{N}$. Then we note that
\begin{equation}
l(s, ~\chi) = 2\sum^{f}_{a=1} \chi(a) H(s, a|f)
\end{equation}
In Eq. (12), we give a value of $l(s, x)$ at negative integer:
\begin{eqnarray*}
l(-n, \chi) &=& 2 \sum_{a=1}^{f} \chi(a)
H(-\eta, ~a|f)\\
&=& 2 \sum_{a=1}^{f} \chi(a) (-1)^{a} \dfrac{f^{n}}{2} E_{n}
\left( \dfrac{a}{f} \right)\\
&=& f^{n} \sum_{a=0}^{f-1} \chi(a) (-1)^{a} E_{n}
\left( \dfrac{a}{f} \right)\\
&=& E_{n, \chi}.
\end{eqnarray*}

The function $H(s, ~a|F)$ will be called partial zeta function
which interpolates Euler polynomials at negative integers. The
values of $l(s, ~\chi)$ at negative integers are algebraic, hence
may be regarded as lying in an extension of $\mathbb{Q}_{p}$. We
therefore look for a $p$-adic function which agrees with $l(s,
~\chi)$ at negative integers in later.
\vspace{0.2in}\\
\noindent{\bf 3. A note on $p$-adic $l$-function}
\vspace{0.1in}\\
\indent We define $\langle x \rangle = \frac{x}{w(x)}$, where
$w(x)$ is the Teichm\"uller character. When $F$(=odd) is a
multiple of $p$ and $(a,p)=1$, we define
$$H_{p}(s, ~a|F) = \dfrac{(-1)^{a}}{2} \langle a \rangle^{-s}
\sum_{j=0}^{\infty} \binom{-s}{j} \left(\dfrac{F}{a}\right)^{j}
E_{j},$$ for $s \in \mathbb{Z}_{p}.$

It is easy to see that
\begin{eqnarray*}
H_{p}(-n, ~a|F) &=& \dfrac{(-1)^{a}}{2} \langle a \rangle^{n}
\sum_{j=0}^{n} \binom{n}{j} \binom{F}{a} E_{j}\\
&=& \dfrac{(-1)^{a}}{2} F^{n} w^{-n}(a) \sum_{j=0}^{n}
\binom{n}{j}\left(\dfrac{a}{F}\right)^{n-j}E_{j}\\
&=& \dfrac{(-1)^{a}}{2} F^{n} w^{-n}(a) E_{n}
\left(\dfrac{a}{F}\right)\\
&=& w^{-n}(a) H(-n, ~a|F),
\end{eqnarray*}
for all positive integers. Now we consider $p$-adic interpolation
function for Euler numbers as follows;
$$l_{p}(s, ~\chi) = 2 \sum_{\substack {a=1 \\ (a, p)=1}}^{F} \chi(a) H_{p}(s, ~a |F)$$
for $s\in \mathbb{Z}_{p}.$

Let $n$ be natural number. Then we have
\begin{eqnarray*}
l_{p}(-n, \chi) &=& 2 \sum_{\substack {n=1 \\ (n, p)=1}}^{F}
\chi(a) H_{p}(-n, a|F)\\
&=& E_{n, ~\chi w^{-n}}-p^{n}\chi w^{-n}(p)E_{n, ~\chi w^{-n}}\\
&=& (1-p^{n} \chi w^{-n}(p)) E_{n, ~\chi w^{-n}}.
\end{eqnarray*}

In fact, we have the formula
$$l_{p}(s, ~\chi) = \sum_{a=1}^{F}(-1)^{a} \langle a \rangle^{-s}
\chi(a) \sum_{j=0}^{\infty} \binom{-s}{j} \left( \dfrac{F}{a}
\right)^{j} E_{j},$$ for $s \in \mathbb{Z}_{p}$.

This is a $p$-adic analytic function and has the following
properties for $\chi = w^{t}$.
$$l_{p}(-n, w^{t}) = (1-p^{n}) E_{n},\quad {\rm where}~~ n \equiv t ~({\rm mod}
~p-1),$$
$$l_{p}(s, w^{t}) \in \mathbb{Z}_{p} ~~{\rm for~ all}~ s \in
\mathbb{Z}_{p},\quad {\rm when}~~ t \equiv 0 ~({\rm mod} ~p-1).$$

If $t \equiv 0 ~({\rm mod}~ p-1)$, then $l_{p}(s_{1}, ~w^{t})
\equiv l_{p}(s_{2}, w^{t}) ~({\rm mod}~p)$ for all $s_{1}, s_{2}
\in \mathbb{Z}_{p}$,
$$l_{p}(k, w^{t}) \equiv l_{p}(k+p, w^{t}) ~({\rm mod}~p).$$

It is easy to see that
$$\dfrac{1}{r+k-1} \binom{-r}{k}
\binom{1-r-k}{j} = \dfrac{-1}{j+k} \binom{-r}{k+j-1}
\binom{k+j}{j},$$ for all positive integers with $r, j, k$ with
$j, k\geq 0$, $j+k > 0$ and $r \neq 1-k$.

Thus, we note that $$\dfrac{1}{r+k-1} \binom{-r}{k}
\binom{1-r-k}{j} = \dfrac{1}{r-1} \binom{-r+1}{k+j}
\binom{k+j}{j}.$$
Hence, we have
$$\dfrac{r}{r+k} \binom{-r+1}{k}
\binom{-r-k}{j} = \binom{-r}{k+j} \binom{k+j}{j},$$ where $k, j$
are positive integers. Let $F$(=odd) be positive integers. Then
\begin{eqnarray*}
&&\sum_{l=0}^{n-1} \dfrac{(-1)^{Fl+a}}{(Fl+a)^{r}} \\
&=&
\sum_{l=0}^{n-1} (-1)^{Fl+a} a^{-r} \sum_{s=0}^{\infty}
\binom{-r}{s} \left( \dfrac{Fl}{a} \right)^{s}\\
&=& \sum_{m=0}^{\infty} \binom{-r}{m} a^{-r}
\left(\dfrac{F}{a}\right)^{m}
(-1)^{a} \sum_{l=0}^{n-1} (-1)^{l} l^{m}\\
&=& \sum_{m=0}^{\infty} \binom{-r}{m} a^{-r} \left( \dfrac{F}{a}
\right)^{m} (-1)^{a} \dfrac{(-1)^{n+1}}{2} \sum_{l=0}^{m-1} E_{l}
n^{m-l} \binom{m}{l} + ((-1)^{m-1} +1)E_{m},
\end{eqnarray*}
when $n$ is even integer
\begin{eqnarray*}
&&\sum_{l=0}^{n-1} \dfrac{(-1)^{Fl+a}}{(Fl+a)^{r}} \\
&=&
a^{-r}(-1)^{a} \sum_{m=0}^{\infty} \binom{-r}{m}
\left(\dfrac{F}{a}\right)^{m} \left\{ \dfrac{(-1)^{m+1}}{2}
\sum_{l=0}^{m+1} E_{l} n^{m-l}
\binom{m}{l} \right\}\\
&=& - a^{-r}(-1)^{a} \sum_{s=0}^{\infty} \binom{-r}{s}
\left(\dfrac{F}{a} \right)^{s} \left(\sum_{l=0}^{s+1} E_{l}
n^{s-l} \binom{s}{l} \right)\\
&=& - \sum_{s=0}^{\infty} \binom{-r}{s} w^{-r}(a)
\left(\dfrac{F}{a} \right)^{s} \dfrac{(-1)^{a}}{2} \langle a
\rangle^{-r} \sum_{l=0}^{s-1} E_{l} n^{s-l} \binom{s}{l}\\
 &=& \sum_{k=0}^{\infty} \sum_{l=0}^{\infty} \binom{-r}{k+l}
w^{-r}(a) \left( \dfrac{a}{F} \right)^{-k-l} n^{k+l}
\dfrac{(-1)^{a}}{2} \langle a \rangle^{-r} E_{l} \left(
\dfrac{1}{n} \right)^{l} \binom{k+l}{l}\\
&=& \sum_{k=0}^{\infty} \sum_{l=0}^{\infty} \dfrac{r}{r+k}
\binom{-r-1}{k} \binom{-r-k}{l} a^{-r} \left( \dfrac{a}{F}
\right)^{-k-l} n^{k+l} \dfrac{(-1)^{a}}{2} E_{l} \left(
\dfrac{1}{n} \right)^{l}\\
&=& \sum_{k=0}^{\infty} \sum_{l=0}^{\infty} \dfrac{r}{r+k}
\binom{-r-1}{k} w^{-k-r}(a) (nF)^{k} \dfrac{(-1)^{a}}{2} \langle a
\rangle^{-k-r} \sum_{l=0}^{\infty} \binom{-r-k}{l}
E_{l}\left(\dfrac{F}{a} \right)\\
&=& - \sum_{k=0}^{\infty} \dfrac{r}{r+k} \binom{-r-1}{k}
w^{-k-r}(a) (nF)^{k} H_{p}(r+k,~a|F).
\end{eqnarray*}

For $F=p,~r,~n$(=even) $\in \mathbb{N}$, we see that
\begin{eqnarray*}
2 \sum_{\substack{j=1 \\ (j, p)=1}}^{np} \dfrac{(-1)^{j}}{j^{p}}
&=& 2 \sum_{a=1}^{p-1} \sum_{l=0}^{n-1}
\dfrac{(-1)^{a+pl}}{(a+pl)^{r}}\\
&=& -2 \sum_{a=1}^{p-1} \sum_{k=0}^{\infty} \dfrac{r}{r+k}
\binom{-r-1}{k} w^{-k-r}(a) (pn)^{k} H_{p}(r+k,~a|p)\\
&=& - \sum_{k=0}^{r} \dfrac{r}{r+k} \binom{-r-1}{k} (pn)^{k} 2
\sum_{a=1}^{p-1}w^{-k-r}(a) H_{p}(r+k,~a|p)\\
&=& -\sum_{k=0}^{r} \binom{-r-1}{k} (pn)^{k} l_{p}(r+k, ~w^{-k-r})\\
&=& -\sum_{k=0}^{\infty} \dfrac{r}{r+k} \binom{-r-1}{k} (pn)^{k}
l_{p}(r+k, ~w^{-k-r}).
\end{eqnarray*}

Therefore we obtain the following:
\vspace{0.1in}\\
\noindent {\bf Theorem 6.} {\it Let $p$ be an odd prime and let
$n($=even$)$, $r \in \mathbb{N}$. Then we have}
$$2 \sum_{\substack{j=1 \\ (j, p)=1}}^{np} \dfrac{(-1)^{j}}{j^{r}}
= - \sum_{k=0}^{\infty} \dfrac{r}{r+k} \binom{-r-1}{k} (pn)^{k}
l_{p}(r+k,~w^{-k-r}).$$
\vspace{0.2in}\\
\footnotesize{

\end{document}